# The Effort of Increasing Reynolds Number in Projection-Based Reduced Order Methods: from Laminar to Turbulent Flows

Saddam Hijazi, Shafqat Ali, Giovanni Stabile,
Francesco Ballarin and Gianluigi Rozza

**Abstract**  We present two different reduced order strategies for incompressible parameterized Navier-Stokes equations characterized by varying Reynolds numbers. The first strategy deals with low Reynolds number (laminar flow) and is based on a stabilized finite element method during the offline stage followed by a Galerkin projection on reduced basis spaces generated by a greedy algorithm. The second methodology is based on a full order finite volume discretization. The latter methodology will be used for flows with moderate to high Reynolds number characterized by turbulent patterns. For the treatment of the mentioned turbulent flows at the reduced order level, a new POD-Galerkin approach is proposed. The new approach takes into consideration the contribution of the eddy viscosity also during the online stage and is based on the use of interpolation. The two methodologies are tested on classic benchmark test cases.

## 1 Introduction

Nowadays we see an increasing need for numerical simulation of fluid dynamics problems with high Reynolds number. These problems come from different types of

Saddam Hijazi
mathLab, Mathematics Area, SISSA, Trieste, Italy, e-mail: `shijazi@sissa.it`

Shafqat Ali
mathLab, Mathematics Area, SISSA, Trieste, Italy, e-mail: `sali@sissa.it`

Giovanni Stabile
mathLab, Mathematics Area, SISSA, Trieste, Italy, e-mail: `gstabile@sissa.it`

Francesco Ballarin
mathLab, Mathematics Area, SISSA, Trieste, Italy, e-mail: `fballarin@sissa.it`

Gianluigi Rozza
mathLab, Mathematics Area, SISSA, Trieste, Italy,e-mail: `grozza@sissa.it`

applications and fields. This pushes the scientific community to offer new techniques and approaches which can meet also the demand of industry to simulate higher Reynolds number fluid problems [11]. Today, in several situations, there is a need to perform simulations in a multi-query contest (e.g. optimization, uncertainty quantification) with an extremely reduced computational time as a requirement (real-time control). Therefore, in such situations, the resolution of the governing PDEs using standard discretization techniques may become unaffordable. Hence, reduced order modelling has become an important tool to reduce the computational complexity. In this chapter, the aim is to present approaches and methodologies to face the problem of building efficient reduced order models (ROMs) for fluid problems with various ranges of the Reynolds number.

This chapter is organized as follows: in section 2 we define the steady Navier-Stokes equations in strong formulation. In section 3 we present residual based stabilized reduced basis method for parameterized Navier-Stokes problem characterized by low Reynolds number. Section 4 deals with POD-Galerkin reduction for parameterized Navier-Stokes problem in case of higher Reynolds number. In section 5 we show some numerical results for both strategies. Finally in section 6 we summarize the main outcomes of this chapter and we outline some perspectives.

## 2 Projection based ROMs

In this section some basic notions of projection based ROMs [13] are recalled. Firstly, the mathematical problem deals with the steady Navier-Stokes equations, and reads as follows:

$$\begin{cases} (\mathbf{u} \cdot \nabla)\mathbf{u} - \nu \Delta \mathbf{u} + \nabla p = \mathbf{0} & \text{in } \Omega, \\ \nabla \cdot \mathbf{u} = 0 & \text{in } \Omega, \\ \mathbf{u} = \mathbf{U}_{in} & \text{on } \partial \Omega_{In}, \\ \mathbf{u} = \mathbf{0} & \text{on } \partial \Omega_0, \\ (\nu \nabla \mathbf{u} - p\mathbf{I})\mathbf{n} = \mathbf{0} & \text{on } \partial \Omega_{Out}, \end{cases} \tag{1}$$

where $\mathbf{u}(\mathbf{x})$ and $p(\mathbf{x})$ are the velocity and pressure fields respectively, $\Omega \subset \mathbb{R}^2$ is a bounded domain, while $\partial \Omega = \partial \Omega_{In} \cup \partial \Omega_0 \cup \partial \Omega_{Out}$ is the boundary of the domain formed by three parts $\partial \Omega_{In}$, $\partial \Omega_0$ and $\partial \Omega_{Out}$ which correspond to the inlet, the physical walls and the outlet respectively, $\mathbf{U}_{in}$ is the velocity at the inlet part of the boundary. $\nu$ is the viscosity of the fluid. Then the problem reads find $\mathbf{u}(\mathbf{x})$ and $p(\mathbf{x})$ which satisfy (1) and lie respectively in the following spaces $\boldsymbol{V} = [\boldsymbol{H}^1(\Omega)]^d$ , and $Q = L_0^2(\Omega)$ see [29] for more details.

In the context of this work, the main goal is studying how the flow fields change as a result of the variation of certain parameters. For this reason a parameterized version of (1) will be considered. The set of parameters is denoted by $\mu$ where this vector of parameters lies in the parameter space $\mathbb{P}$, note that $\mathbb{P}$ is compact set in $\mathbb{R}^p$ with $p$ being the length of the vector $\mu$. The parameters can be geometrical or physical parameters or a combination of them [7]. The objective is to be able to compute the velocity and pressure fields for every parameter value inside the parameter space. The cost of doing that operation resorting on full order methods can be prohibitive. For this reason ROMs [12, 30, 7, 18] have been developed, as an approach to achieve the objective of computing efficiently and accurately the flow fields, when the input parameters are being varied.

One key assumption in ROMs is that the dynamics of the system under study is governed by a reduced number of dominant modes. In other words, the solution to the full order problem lies in a low dimensional manifold that is spanned by the previously mentioned modes [18]. Consequently the velocity and pressure fields can be approximated by decomposing them into linear combination of global basis functions $\phi_i(\mathbf{x})$ and $\chi_i(\mathbf{x})$ (which do not depend on $\mu$) multiplied by unknown coefficients $a_i(\mu)$ and $b_i(\mu)$, for velocity and pressure respectively, then this approximation reads as follows:

$$\mathbf{u}(\mathbf{x};\mu) \approx \sum_{i=1}^{N_u} a_i(\mu)\phi_i(\mathbf{x}), \quad p(\mathbf{x};\mu) \approx \sum_{i=1}^{N_p} b_i(\mu)\chi_i(\mathbf{x}). \tag{2}$$

The reduced basis spaces $\mathbb{V}_{rb} = \mathrm{span}\,\{\phi_i\}_{i=1}^{N_u}$ and $Q_{rb} = \mathrm{span}\,\{\chi_i\}_{i=1}^{N_p}$ can be obtained either by Reduced Basis (RB) method with a greedy approach [18], using Proper Orthogonal Decomposition (POD) [36], by The Proper Generalized Decomposition [15], or by Dynamic Mode Decomposition [33]. In the next two sections we will consider RB and POD methods.

## 3 Stabilized Finite Element RB Reduced Order Method

In this section, we present a RB method for parameterized steady Navier-Stokes problem [30] which ensures stable solution [2]. Our focus in this section is to deal with flows at low Reynolds number with particular emphasis on inf-sup stability at reduced order level.

We know that the Galerkin projection on RB spaces does not guarantee the fulfillment of equivalent reduced inf-sup condition [31]. To fulfill this condition we have to enrich the RB velocity space with the solutions of a supremizer problem [32, 5]. In this work we propose a residual based stabilization technique which circumvents the inf-sup condition and guarantees stable RB solution. This approach consists in adding some stabilization terms into the Galerkin finite element formulation of (1) using equal order ($\mathbb{P}_k/\mathbb{P}_k; k = 1, 2$) velocity pressure interpolation, and than projecting onto RB spaces. As the results in section 5 will show, residual based stabilization methods improves the stability of Galerkin finite element method without compromising the consistency.

We start with introducing two finite-dimensional subspaces $\mathbf{V}_h \subset \mathbf{V}$, $Q_h \subset Q$ of dimension $\mathcal{N}_u$ and $\mathcal{N}_p$, respectively, being $h$ related to the computational mesh size.

The Galerkin finite element approximation of the parameterized problem (1) with the addition of stabilization terms reads as follows: for a given parameter value $\mu \in \mathbb{P}$, we look for the full order solution $(\mathbf{u}_h(\mu), p_h(\mu)) \in \mathbf{V}_h \times Q_h$ such that

$$
\begin{cases}
a(\mathbf{u}_h, \mathbf{v}_h; \mu) + c(\mathbf{u}_h, \mathbf{u}_h, \mathbf{v}_h; \mu) + b(\mathbf{v}_h, p_h; \mu) = \xi_h(\mathbf{v}_h; \mu) & \forall \mathbf{v}_h \in \mathbf{V}_h, \\
b(\mathbf{u}_h, q_h; \mu) = \psi_h(q_h; \mu) & \forall q_h \in Q_h,
\end{cases}
\tag{3}
$$

which we name as the stabilized Galerkin finite element formulation, where $a(.,.;\mu)$ and $b(.,.;\mu)$ are the bilinear forms related to diffusion and pressure-divergence operators, respectively and $c(.,.,.;\mu)$ is the trilinear form related to the convective term. The stabilization terms $\xi_h(\mathbf{v}_h; \mu)$ and $\psi_h(q_h; \mu)$ are defined as:

$$
\begin{aligned}
\xi_h(\mathbf{v}_h; \mu) &:= \delta \sum_K h_K^2 \int_K (-\nu \Delta \mathbf{u}_h + \mathbf{u}_h \cdot \nabla \mathbf{u}_h + \nabla p_h, -\gamma \nu \Delta \mathbf{v}_h + \mathbf{u}_h \cdot \nabla \mathbf{v}_h), \\
\psi_h(q_h; \mu) &:= \delta \sum_K h_K^2 \int_K (-\nu \Delta \mathbf{u}_h + \mathbf{u}_h \cdot \nabla \mathbf{u}_h + \nabla p_h, \nabla q_h),
\end{aligned}
\tag{4}
$$

where $K$ is an element of the domain, $h_K$ is the diameter of $K$, $\delta$ is the stabilization coefficient such that, $0 < \delta \le C$ (C is a suitable constant) needs to be chosen properly [24, 9]. For $\gamma = 0, 1, -1$, the stabilization (4) is respectively known as Streamline Upwind Petrov Galerkin (SUPG) [10], Galerkin least-squares (GLS) [20] and Douglas-Wang (DW) [14].

Next step is to construct the RB spaces $\mathbf{V}_{rb}$ and $Q_{rb}$, for velocity and pressure, respectively. These spaces are constructed using the greedy algorithm [18] and may or may not be enriched with supremizer [2]. In order to control the condition number of RB matrix, the basis functions $\phi_i(\mathbf{x})$ and $\chi_i(\mathbf{x})$ for RB velocity and pressure, respectively are orthonormalized by using the Gram-Schmidt orthonormalization process [18].

Now we write the RB formulation, i.e, we perform a Galerkin projection of (3) onto the RB spaces. Therefore the reduced problem reads as follows: for any $\mu \in \mathbb{P}$, find $(\mathbf{u}_N(\mu), p_N(\mu)) \in \mathbf{V}_{rb} \times Q_{rb}$ such that

$$
\begin{cases}
a(\mathbf{u}_N, \mathbf{v}_N; \mu) + c(\mathbf{u}_N, \mathbf{u}_N, \mathbf{v}_N; \mu) + b(\mathbf{v}_N, p_N; \mu) = \xi_N(\mathbf{v}_N; \mu) & \forall \mathbf{v}_N \in \mathbf{V}_{rb}, \\
b(\mathbf{u}_N, q_N; \mu) = \psi_N(q_N; \mu) & \forall q_N \in Q_{rb},
\end{cases}
\tag{5}
$$

where $\xi_N(\mathbf{v}_N; \mu)$ and $\psi_N(q_N; \mu)$ are the reduced order counterparts of the stabilization terms defined in (4). We call (5) as the stabilized RB formulation.

The Galerkin projection of (3) onto RB spaces can also be performed without adding the stabilization terms in RB formulation. Therefore we have two options here [28]; the first option is the *offline-online stabilization*, where we apply the Galerkin projection on stabilized formulations in both the offline and the online stages, and the second option is *offline-only stabilization*, where we apply stabilization only in the offline stage and then we perform the online stage using the standard formulation. Finally, combining these two options with the supremizer enrichment [32], we come up with the following four options to discuss [2]:

- *offline-online stabilization* with supremizer;
- *offline-online stabilization* without supremizer;
- *offline-only stabilization* with supremizer;
- *offline-only stabilization* without supremizer.

An extension of the work presented in this section to unsteady problems is currently in progress [3].

## 4 Finite Volume POD-Galerkin reduced order model

In this section, the treatment of flow with high Reynolds number will be addressed. The starting point is with the POD-Galerkin projection method in the first subsection, and then the ROM for turbulent flows will be proposed in the second subsection.

### 4.1 POD-Galerkin Projection Method

POD is a very popular method for generating reduced order spaces. It is based on constructing a reduced order space which is optimal in the sense that it minimizes the projection error (the $L^2$ norm of the difference between the snapshots and their projection onto the reduced order basis). After generating the POD space one can project (1) into that space. This approach is called POD-Galerkin projection which has been widely used for building ROMs for variety of problems in Computational Fluid Dynamics (CFD) [27, 21, 8, 1, 16, 4].

The POD space is obtained by solving the following minimization problem :

$$\mathbb{V}_{POD} = \arg\min \frac{1}{N_s} \sum_{n=1}^{N_s} ||\mathbf{u}_n - \sum_{n=1}^{N_s} (\mathbf{u}_n, \phi_i)_{L^2(\Omega)} \phi_i||^2_{L^2(\Omega)}, \quad (6)$$

where $\mathbf{u}_n$ is a general snapshot of the velocity field which is obtained for the sample $\mu_n$ and $N_s$ is the total number of snapshots. The minimization problem can be solved by performing Singular Value Decomposition (SVD) on the matrix formed by the snapshots, or by computing a correlation matrix whose entries are the scalar product between the snapshots and then performing eigenvalue decomposition on that correlation matrix, for more details we refer the reader to [34, 23].

The next step in building the reduced order model is to project the momentum equation of (1) onto the POD space spanned by the velocity POD modes, namely:

$$(\phi_i, (\mathbf{u} \cdot \nabla)\mathbf{u} - \nu \Delta \mathbf{u} + \nabla p)_{L^2(\Omega)} = 0. \quad (7)$$

Inserting the approximations (2) into (7) yields the following system:

$$\nu \mathsf{B} \mathbf{a} - \mathbf{a}^T \mathsf{C} \mathbf{a} - \mathsf{H} \mathbf{b} = \mathbf{0}, \tag{8}$$

where $\mathbf{a}$ and $\mathbf{b}$ are the vectors of coefficients $a_i(\mu)$ and $b_i(\mu)$ , respectively, while the other terms are computed as follows :

$$B_{ij} = (\phi_i, \Delta \phi_j)_{L^2(\Omega)}, \tag{9}$$

$$C_{ijk} = (\phi_i, \nabla \cdot (\phi_j \otimes \phi_k)))_{L^2(\Omega)}, \tag{10}$$

$$H_{ij} = (\phi_i, \nabla \chi_j)_{L^2(\Omega)}. \tag{11}$$

For better understanding of the treatment of the nonlinearity introduced by the convective term the reader may refer to [34]. To close the system (8) an additional number of $N_p$ equations is needed since there are just $N_u$ equations but with $N_u + N_p$ unknowns. The continuity equation cannot be directly used to close the system since the snapshots which are obtained using the full order solver are already divergence free, and the velocity POD modes which are obtained using those snapshots have the same property. This problem can be overcome by two possible approaches, the first one is to use Poisson equation for pressure to get the needed additional equations such that one can close the system. Poisson equation for pressure can be derived by just taking the divergence of the momentum equation and then exploiting the continuity equation. The second possible approach, is the supremizer stabilization method [5, 32] which has been already mentioned in section 3. The latter approach has been developed for finite volume discretization method as well and one can refer to [36] for more details on that. The supremizer approach will ensure that the velocity modes are not all divergence free. One can project the continuity equation onto the space spanned by the pressure modes, which results in the following system:

$$\begin{cases} \nu \mathsf{B} \mathbf{a} - \mathbf{a}^T \mathsf{C} \mathbf{a} - \mathsf{H} \mathbf{b} = \mathbf{0}, \\ \mathsf{P} \mathbf{a} = \mathbf{0}, \end{cases} \tag{12}$$

where the new matrix is $\mathsf{P}$, is computed as follows :

$$P_{ij} = (\chi_i, \nabla \cdot \phi_j)_{L^2(\Omega)}. \tag{13}$$

Concerning the treatment of boundary conditions, a lifting function method is employed. A new set of snapshots with homogeneous boundary condition is created. For the selection of an appropriate lifting function, several options are available such as snapshots average or the solution to a linear problem. We decided here to rely on the latter approach. For more details one can refer to [34].

### 4.2 POD-Galerkin Reduced Order Model for Turbulent Flows

In this subsection, the main goal is to focus on flows which have higher Reynolds number than those considered in the section 3. In these flows the turbulence phenomenon is present. The full order discretization technique used in this case for solving (1) is the Finite Volume Method (FVM) [26, 37] which is widely used in industrial applications. One advantage of the FVM is that the equations are written in conservative form, and therefore the conservation law is ensured at a local level.

The turbulence modelling is employed using $k - \omega$ turbulence model [25] which is a two equations model, is used to ensure the stability of the simulation. In this model the eddy viscosity $v_t$ depends algebraically on two variables $k$ and $\omega$ respectively stand for the turbulent kinetic energy and the specific turbulent dissipation rate. The values of these two variables are computed solving two additional PDEs. The new set of equations to be solved is the Reynolds Averaged Navier-Stokes (RANS) equations which read as follow:

$$\begin{cases} (\boldsymbol{u} \cdot \nabla) \boldsymbol{u} = \nabla \cdot \left[ -p\mathbf{I} + (v + v_t) \left( \nabla \boldsymbol{u} + (\nabla \boldsymbol{u})^T \right) - \frac{2}{3} k \mathbf{I} \right], \\ \nabla \cdot \boldsymbol{u} = 0, \\ v_t = f(k, \boldsymbol{\omega}), \\ \text{Transport-Diffusion equation for } k, \\ \text{Transport-Diffusion equation for } \omega. \end{cases} \quad \text{in } \Omega$$

In order to build a reduced order model for the new set of equations one can extend the previous assumption (2) to the eddy viscosity field, namely:

$$v_t(\mathbf{x}; \mu) \approx \sum_{i=1}^{N_{v_t}} g_i(\mu) \eta_i(\mathbf{x}),$$

The eddy viscosity modes $\eta_i$ are computed similarly to those of velocity and pressure. Following the procedure explained in subsection 4.1 one can project the momentum equation onto the spatial bases of velocity. The continuity equation is projected onto the spatial bases of pressure with the use of a supremizer stabilization approach. In contrast, $k - \omega$ transport-diffusion equations are not used in the projection procedure, this makes the reduced order model general and independent of the turbulence model used in the full order simulations.

The resulting system is the following:

$$\begin{cases} (\mathsf{B} + \mathsf{B}_\mathsf{T})\mathbf{a} - \mathbf{a}^\mathbf{T} \mathsf{C} \mathbf{a} + \mathbf{g}^\mathbf{T} (\mathsf{C}_{\mathsf{T}1} + \mathsf{C}_{\mathsf{T}2})\mathbf{a} - \mathsf{H}\mathbf{b} = \mathbf{0}, \\ \mathsf{P}\mathbf{a} = 0, \end{cases} \quad (14)$$

Where the new terms with respect to the dynamical system in (12) are computed as follows:

$$B_{T_{ij}} = \left( \phi_i, \nabla \cdot (\nabla \phi_j^T) \right)_{L^2(\Omega)}, \tag{15}$$

$$C_{T1_{ijk}} = \left( \phi_i, \eta_j \Delta \phi_k \right)_{L^2(\Omega)}, \tag{16}$$

$$C_{T2_{ijk}} = \left( \phi_i, \nabla \cdot \eta_j (\nabla \phi_k^T) \right)_{L^2(\Omega)}. \tag{17}$$

One can see that a new set of coefficients $g_i$ has been introduced. These coefficients are used in the approximation of the eddy viscosity fields, and in order to compute them an interpolation procedure using Radial Basis Functions (RBF) [22] has been used in the online stage. After that one can solve the system (14) for the vectors of coefficients **a** and **b**.

In the remaining part of this subsection, the interpolation method used to compute the coefficients of the reduced viscosity will be explained in further details. The starting point consists of the set of samples used in the offline stage $X_\mu = \{\mu_1, \mu_2, ..., \mu_{N_s}\}$. The associated outputs $y_i$ are the coefficients resulted from the projection of the viscosity snapshots that correspond to each $\mu_i$ onto the viscosity spatial modes $[\chi_j]_{j=1}^{N_{v_t}}$. The goal is to interpolate the known coefficients by making the use of RBF $\zeta_i$ for $i = 1, ..., N_s$. One may assume that $Y$ has the following form:

$$Y(x) = \sum_{j=1}^{N_s} w_j \zeta_j (\|x - x_j\|_2), \tag{18}$$

where $w_j$ are some appropriate weights. In order to interpolate the known data, the following property is required:

$$Y(x_i) = y_i, \quad \text{for} \quad i = 1, 2, ..., N_s. \tag{19}$$

In other words,

$$\sum_{j=1}^{N_s} w_j \zeta_j (\|x_i - x_j\|_2) = y_i, \quad \text{for} \quad i = 1, 2, ..., N_s. \tag{20}$$

The latter system can be solved to find the weights. The procedure dealing with what concerns the use of RBF interpolation is summarized in the following box for both offline and online stages.

In the context of this work, RBF is used according to the following algorithm. The methodology has two parts, the first is within the offline stage in which the interpolant RBF is constructed. The second part, which takes place during the online stage, consists into the evaluation of the coefficients $[g_i]_{i=1}^{N_{v_t}}$ using the latter mentioned RBF methodology.

**Offline Stage**

**Input**: The set of samples for which the offline stage has been run $X_\mu = \{\mu_1, \mu_2, ..., \mu_{N_s}\}$, with the corresponding eddy viscosity snapshots $\nu_{t1}, \nu_{t2}, ..., \nu_{tN_s}$, the number of eddy viscosity modes to be used in the re-duction during online phase $N_{v_t}$ and finally $i = 1$ which is an index to be used during the stage.

**Goal**: for $i = 1, 2, ..., N_{v_t}$ construct $g_i(\mu) = \sum_{j=1}^{N_s} w_{i,j} \zeta_{i,j}(\|\mu - \mu_j\|_2)$

**Step 1**

Compute the eddy viscosity modes $[\chi_k]_{k=1}^{N_{v_t}}$ using POD as mentioned before.

**Step 2**

Compute the following coefficients

$$g_{i,j} = (v_{t\,j}, \chi_i)_{L^2(\Omega)}, \quad \text{for} \quad j = 1, 2, ... N_s \quad . \tag{21}$$

**Step 3**

Solve the following linear system for the vector of weights $\mathbf{w}_i = [w_{i,j}]_{j=1}^{N_s}$

$$\sum_{j=1}^{N_s} w_{i,j} \zeta_{i,j}(\|\mu_k - \mu_j\|_2) = g_{i,k}, \quad \text{for} \quad k = 1, 2, ..., N_s. \tag{22}$$

**Step 4**

Store the weights $[w_{i,j}]_{j=1}^{N_s}$ and construct the scalar coefficients $g_i(\mu)$.

**Step 5**

If $i = N_{v_t}$ **terminate**, otherwise set $i = i + 1$ and go to **Step 2**.

**Online Stage**

As **Input** we have the new parameter value $\mu^*$ and the goal is to compute $\mathbf{g}(\mu^*) = [g_i(\mu^*)]_{i=1}^{N_{v_t}}$

Which is done simply by computing $g_i(\mu^*) = \sum_{j=1}^{N_s} w_{i,j} \zeta_{i,j}(\|\mu^* - \mu_j\|_2)$ for $i = 1, 2, ..., N_{v_t}$

After computing the coefficients of the viscosity reduced order solution $[g_i]_{i=1}^{N_{v_t}}$ then it will be possible to solve the reduced order system (14). Afterwards, one can compute the reduced order solution for both velocity and pressure using (2). From now on this approach will be referred to as *POD-Galerkin-RBF ROM*. The POD-Galerkin-RBF model will be tested on a simple benchmark test case of the backstep in steady setting, with the offline phase being done with a RANS approach. For the application of this model on more complex cases involving LES full order simulations and in an unsteady setting the reader may refer to [19].

## 5 Numerical Results

In this section we present numerical results for both reduced order modelling strategies presented in the previous sections. In subsection 5.1 we present the numerical results for low Reynolds number using stabilized RB method developed in section 3 for steady Navier-Stokes equations. Subsection 5.2 is based on the results for POD-Galerkin-RBF on a backward facing step problem. In both cases we consider only physical parameters.

### 5.1 Stabilized Finite Element based ROM Results

In this test case, we apply the stabilized RB model developed in section 3 for the Navier-Stokes problem to the *lid driven-cavity* problem with only one physical parameter $\mu$ which denotes the Reynolds number. We consider only the first three options and we have done several test cases to compare the three options. Fourth option is the worst option and is not reported here. The stabilization option that we consider here is the SUPG stabilization, corresponding to $\gamma = 0$ in (4). The computational domain is shown in Fig. 1 and the boundary conditions are

$$u_1 = 1, u_2 = 0 \text{ on } \partial \Omega_{In} \text{ and } \mathbf{u} = \mathbf{0} \text{ on } \partial \Omega_0 \tag{23}$$

The mesh of this problem is non-uniform with 3794 triangles and 1978 nodes,

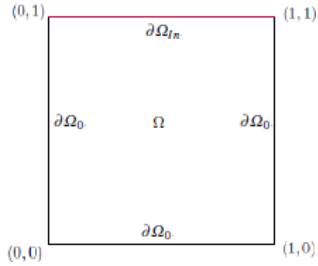

**Fig. 1** Unit cavity domain $\Omega$ for RB problem with boundaries identified.

whereas the minimum and maximum size elements are $h_{min} = 0.0193145$ and $h_{max} = 0.0420876$, respectively. All the numerical simulations for this case are performed using FreeFem++ [17] and RBniCS [6].

Figure 2 shows the FE velocity (left), RB velocity obtained using *offline-online stabilization* (center), and the RB velocity obtained for *offline-only stabilization* (right). From these solutions we see that the FE and RB solutions are similar.

Figure 3 plots the FE pressure (left), RB pressure obtained using *offline-online stabilization* (center), and the RB pressure obtained for *offline-only stabilization* (right). These results show that the RB pressure with *offline-online stabilization* is stable but RB solution obtained by *offline-only stabilization* is highly oscillatory even with the supremizer enrichment. All these solutions are obtained for equal order linear velocity pressure interpolation $\mathbb{P}_1/\mathbb{P}_1$. Similar results can be shown for $\mathbb{P}_2/\mathbb{P}_2$ [2].

Figure 4 illustrates the error between FE and RB solutions for velocity (left) and pressure (right). We show the comparison between *offline-online stabilization* with/without supremizer and *offline-only stabilization* with supremizer. These comparison shows that the *offline-online stabilization* is the most appropriate way to stabilize and the enrichment of RB velocity space with supremizer may not be necessary. We are getting even a better approximation of the velocity without the supremizer, which is polluted a little bit by the supremizer. However in case of pressure, supremizer is improving the accuracy in the case of *offline-online stabilization*. All the results here are presented for equal order linear velocity pressure interpolation $\mathbb{P}_1/\mathbb{P}_1$.

In Table 1 we summarize the computational cost of offline and online stage for different choices of FE spaces, parameter detail, FE and RB dimensions. From this table we can see that the *offline-online stabilization* without supremizer is less expensive as compared to *offline-online stabilization* with supremizer.

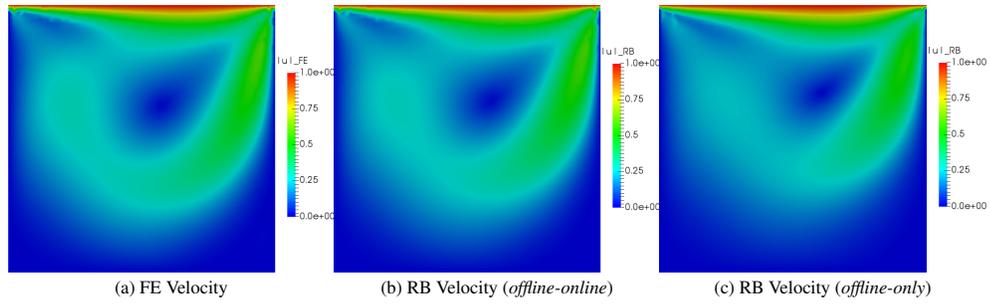

**Fig. 2** SUPG stabilization: FE and RB solutions for velocity at $Re = 200$.

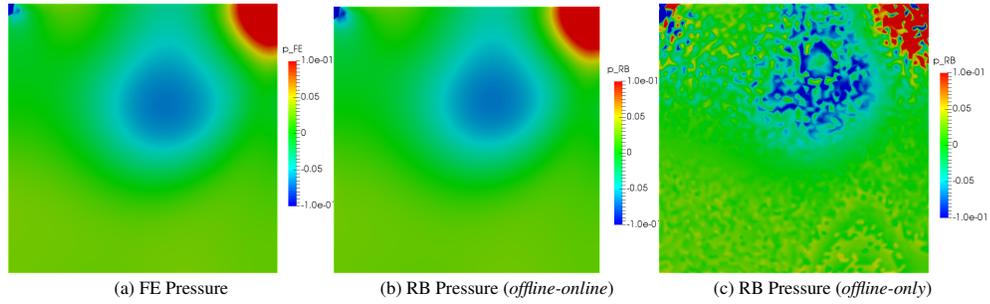

**Fig. 3** SUPG stabilization: FE and RB solutions for pressure at $Re = 200$.

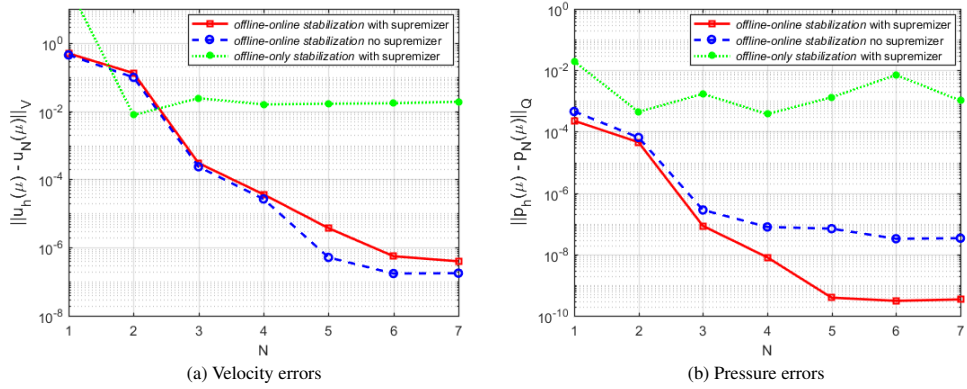

(a) Velocity errors                    (b) Pressure errors

**Fig. 4** Error between FE and RB solutions: velocity (left) and pressure (right), obtained by different options using SUPG stabilization.

**Table 1** Computational details of steady Navier-Stokes problem with physical parameter only.

| Physical parameter | $\mu$ (Reynolds number) |
|---|---|
| Range of $\mu$ | [100,500] |
| *Online* $\mu$ | 200 |
| FE degrees of freedom | 13218 ($\mathbb{P}_1/\mathbb{P}_1$) |
| | 52143 ($\mathbb{P}_2/\mathbb{P}_2$) |
| RB dimension | $N_u = N_s = N_p = 7$ |
| *Offline* time ($\mathbb{P}_1/\mathbb{P}_1$) | 1182$s$ (*offline-online stabilization* with supremizer) |
| | 842$s$ (*offline-online stabilization* without supremizer) |
| *Offline* time ($\mathbb{P}_2/\mathbb{P}_2$) | 2387$s$ (*offline-online stabilization* with supremizer) |
| | 2121$s$ (*offline-online stabilization* without supremizer) |
| *Online* time ($\mathbb{P}_1/\mathbb{P}_1$) | 74$s$ (with supremizer) |
| | 65$s$ (without supremizer) |
| *Online* time ($\mathbb{P}_2/\mathbb{P}_2$) | 131$s$ (with supremizer) |
| | 108$s$ (without supremizer) |

## 5.2 Finite Volume POD-Galerkin-RBF ROM Results

In this subsection the numerical results for the reduced order model obtained using the POD-Galerkin-RBF approach are shown. The finite volume C++ library OpenFOAM® (OF) [38] is used as the numerical solver at the full order level. At the reduced order level the reduction is done using the library ITHACA-FV [35] which is based on C++.

We have tested the proposed model on the benchmark case of the backstep see Fig. 5. The test is performed in steady state setting, the two considered parameters are both physical and consist into the magnitude of the velocity at the inlet and the inclination of the velocity with respect to the inlet. In addition a comparison is presented between the results obtained using the newly developed POD-Galerkin-RBF approach with the POD-Galerkin option that is not using RBF in the online stage.

The interest of this test is in reducing the Navier-Stokes equations in the case of turbulent flows or flows with high Reynolds number. In this case the value of the Reynolds number is around $10^4$, while the physical viscosity $\nu$ is equal to $10^{-3}$. $\mu = [\mu_1, \mu_2]$ is the vector of the parameters with $\mu_1$ being the magnitude of the velocity at the inlet and $\mu_2$ the inclination of the velocity with respect to the inlet which is measured in degrees. Samples for both parameters are generated as 20 equally distributed points in the ranges of $[0.18, 0.3]$ and $[0, 30]$ respectively. The training of the reduced order model is done with the generated 400 sample points in the offline stage. The Reynolds number as mentioned before is of order of $10^4$ and ranges from $9.144 \times 10^3$–$1.524 \times 10^4$.

In Fig. 5 one can see the computational domain that has been used in this work. The characteristic length $d$ and is equal to 50.8 meters. In the full order problem the boundary conditions for velocity and pressure are set as reported in Table 2.

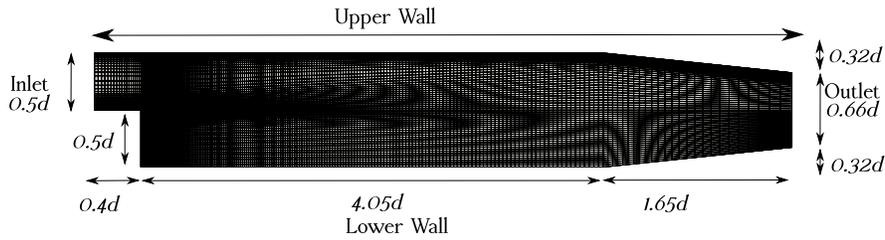

**Fig. 5** The computational domain used in the numerical simulations, $d$ is equal to 50.8 meters.

In the reduced order model the supremizer approach has been used to stabilize pressure. In Table 3 one can see the cumulative eigenvalues for velocity, pressure, supremizer (which is denoted by **S**) and viscosity.

During the online phase another set of samples has been used to check the reduced order model which is a cross validation test of the model. The value of the parameter vector given in the online phase is denoted by $\mu_i^*$ where $i = 1, ..., N_{online-samples}$. The samples which were used in the cross validation have been chosen inside the ranges of the samples used in the offline stage. For the sake of better evaluation of the model the online samples have been chosen such that they are as far as possible in the parameter space from those used to educate the model. After taking that criterion into consideration the samples used in the online phase happened to be equally distributed in the ranges of $[0.20826, 0.28405]$ and $[5.5368, 29.221]$ for $\mu_1$ and $\mu_2$, respectively. Seven samples were used for $\mu_1$ and six for $\mu_2$.

Recall that in this case the parameters were basically the two components of the velocity at the inlet. Therefore two lifting fields were computed which correspond to the full order solution for the velocity field with unitary boundary condition. The first and second lifting fields are the steady state solutions with the velocity at the inlet being $U = (1,0)$ and $U = (0,1)$, respectively. These two fields are added to the velocity modes.

The RBF functions for the turbulent viscosity are chosen to be Gaussian functions. The system (14) has been solved for each online sample $\mu_i^*$ in the online phase and the fields have been constructed. The ROM fields obtained by solving (14) have been compared to those resulted from solving the POD-Galerkin system (12), which does not take into consideration the contribution of the eddy viscosity.

In Fig. 6 one can see the velocity fields obtained by the full order solver, ROM velocity field obtained by the POD-Galerkin approach and ROM velocity obtained by the new POD-Galerkin-RBF model. In Fig. 7 there is the same comparison but for the pressure fields. In both figures the online sample which has been introduced to both reduced order models is the one with $\mu^* = (0.22089, 24.484)$, which corresponds to the velocity vector at the inlet to be $\mathbf{U} = (0.20103, 0.091548)$. The reduction has been made with seven modes for velocity, pressure, supremizer and eddy viscosity (just considered in the POD-Galerkin-RBF model). One can see that the POD-Galerkin-RBF model is able to capture the dynamics efficiently. It has successfully reconstructed the full order solution from both qualitative and quantitative aspects. On the other hand, it is quite clear that the classical POD-Galerkin model, which does not consider the contribution of the eddy viscosity in its formulation fail to give an accurate reproduction of the full order solution, especially close to the top and to the outlet for the velocity field and at the inlet for the pressure field.

Looking on the results from a quantitative point of view, in the POD-Galerkin-RBF model we have values of 0.00612 and 0.02957 for the relative error in $L^2$ norm for velocity and pressure, respectively, while the POD-Galerkin model has errors of 0.37967 and 2.2296. Table 4 shows a comparison between the two models in terms of the error over all the samples used in the online phase (average and maximum value). Figures 8 and 9 show the error as function of the two parameters when one of them is fixed and the other is varied.

**Table 2**  Boundary Conditions

|   | inlet | outlet | lower and upper walls |
|---|---|---|---|
| $\boldsymbol{u}$ | $\mathbf{u_{in}} = [\mu_1 cos(\mu_2), \mu_1 sin(\mu_2)]$ | $\nabla\mathbf{u}\cdot\mathbf{n} = \mathbf{0}$ | $\mathbf{u} = \mathbf{0}$ |
| $p$ | $\nabla p \cdot \mathbf{n} = 0$ | $p = 0$ | $\nabla p \cdot \mathbf{n} = 0$ |

# 6 Conclusion and Perspectives

In this chapter we have proposed two different ROM strategies for the incompressible parameterized Navier-Stokes equations to deal from low to higher Reynolds

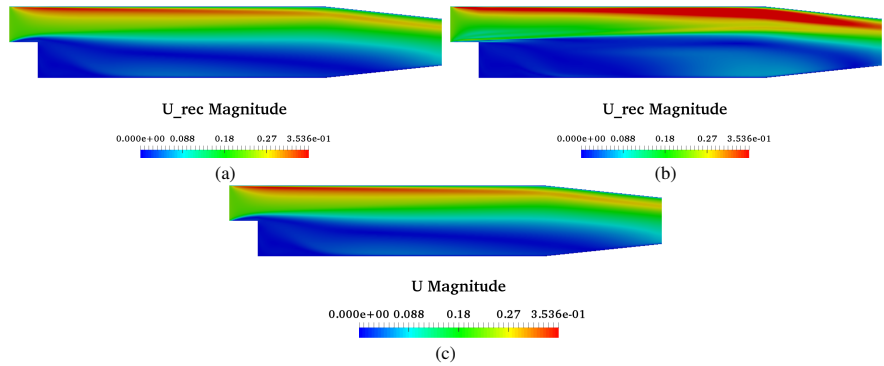

**Fig. 6** Velocity fields: (a) shows the ROM Velocity obtained by POD-Galerkin-RBF ROM model, while in (b) one can see the ROM Velocity (without viscosity incorporated in ROM), and finally in (c) we have the FOM Velocity.

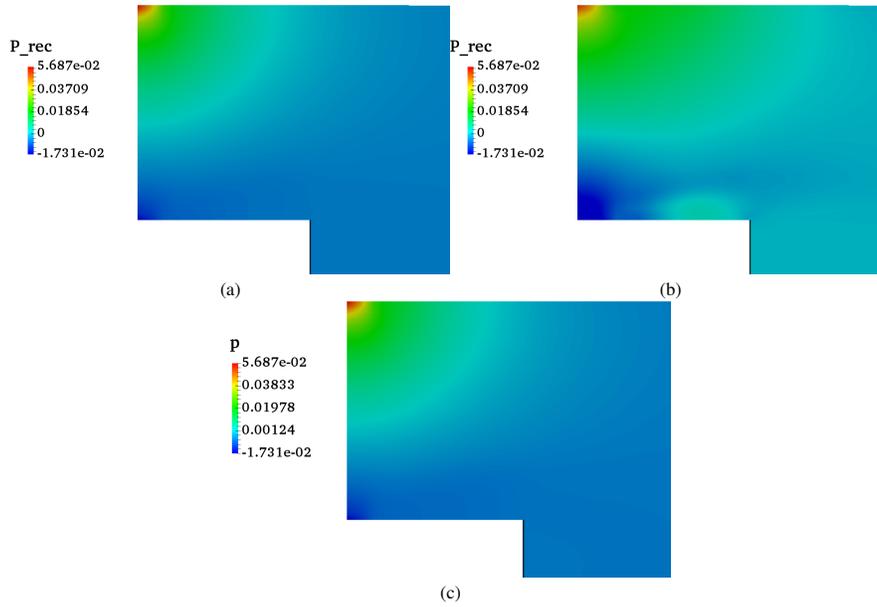

**Fig. 7** Pressure fields: (a) shows the ROM Pressure obtained by POD-Galerkin-RBF ROM model, while in (b) one can see the ROM Pressure (without viscosity incorporated in ROM), and finally in (c) we have the FOM Pressure.

16       Saddam Hijazi, Shafqat Ali, Giovanni Stabile, Francesco Ballarin and Gianluigi Rozza

**Table 3** Cumulative Eigenvalues

| N Modes | $\mathbf{u}$ | $p$ | $\mathbf{S}$ | $\nu_t$ |
|---|---|---|---|---|
| 1 | 0.971992 | 0.868263 | 0.899488 | 0.985703 |
| 2 | 0.993017 | 0.998541 | 0.996392 | 0.998884 |
| 3 | 0.997589 | 0.999915 | 0.999767 | 0.999673 |
| 4 | 0.999196 | 0.999963 | 0.999929 | 0.999880 |
| 5 | 0.999545 | 0.999985 | 0.999965 | 0.999926 |
| 6 | 0.999828 | 0.999997 | 0.999988 | 0.999971 |
| 7 | 0.999914 | 0.999999 | 0.999996 | 0.999986 |
| 8 | 0.999952 | 0.999999 | 0.999998 | 0.999992 |
| 9 | 0.999978 | 1.000000 | 0.999999 | 0.999995 |
| 10 | 0.999986 | 1.000000 | 0.999999 | 0.999997 |

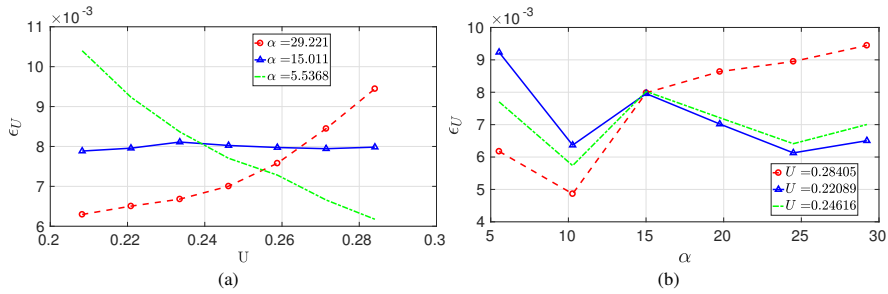

**Fig. 8** The $L^2$ relative error for velocity fields as function of the parameters. In (a) the error is plotted versus the inclination of the velocity at the inlet. While in (b) the error is plotted versus the magnitude of the velocity at the inlet.

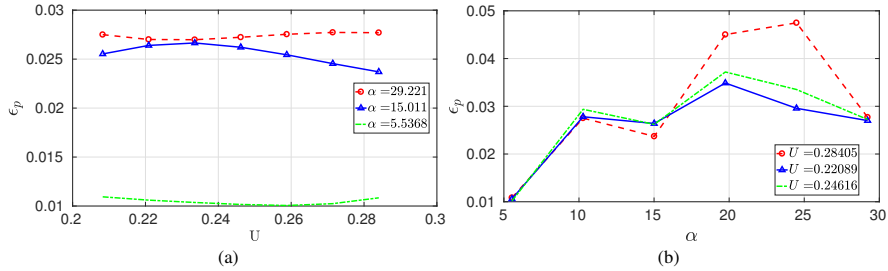

**Fig. 9** The $L^2$ relative error for pressure fields as function of the parameters. In (a) the error is plotted versus the inclination of the velocity at the inlet. While in (b) the error is plotted versus the magnitude of the velocity at the inlet.

**Table 4** Relative $L^2$ Error for Velocity and Pressure Fields: *Average* is taken over all samples used in the online phase, while *maximum* represents the worse case among the samples. POD-Galerkin-RBF model results are compared to those of the normal POD-Galerkin one.

| | $\mathbf{u}$ with RBF | $p$ with RBF | $\mathbf{u}$ without RBF | $p$ without RBF |
|---|---|---|---|---|
| Average Relative Error | 0.0073 | 0.0276 | 0.2592 | 1.5412 |
| Maximum Relative Error | 0.0104 | 0.0475 | 0.3810 | 2.3616 |

number, respectively. In case of low Reynolds number, we have used a stabilized FE discretization techniques at the full order level and then we performed Galerkin projection onto RB spaces, obtained by a greedy algorithm. We have compared the *offline-online stabilization* approach with supremizer enrichment in context of RB inf-sup stability. Based on numerical results, we conclude that a residual based stabilization technique, if applied in both offline and online stage (*offline-online stabilization*), is sufficient to ensure a stable RB solution and therefore we can avoid the supremizer enrichment which consequently reduces the online computation cost. Supremizer may help in improving the accuracy of pressure approximation. We also conclude that a stable RB solution is not guaranteed if we stabilize the offline stage and not the online stage (*offline-only stabilization*) even with supremizer enrichment.

For higher Reynolds number, the test case was the backstep benchmark test case, we have used the FV discretization technique at the full order level. At the reduced order level, we have used a POD-Galerkin projection approach taking into consideration the contribution of the eddy viscosity. The newly proposed approach involves the usage of radial basis functions interpolation in the online stage. The model has been tested on the benchmark case of the backstep, the results showed that the proposed model has successfully reduced RANS equations. On the other hand, the classical POD-Galerkin approach has not been able to reduce the equations accurately in the same study case.

For the future work, we aim to extend the POD-Galerkin-RBF approach to work also on unsteady Navier-Stokes equations. In addition, one important goal is to reduce problems where the offline phase is simulated with LES.

`acknowledgement` We acknowledge the support provided by the European Research Council Consolidator Grant project Advanced Reduced Order Methods with Applications in Computational Fluid Dynamics - GA 681447, H2020-ERC COG 2015 AROMA-CFD, and INdAM-GNCS.